\documentclass{amsart}
\numberwithin{equation}{section}
\usepackage{amssymb}
\usepackage{amsmath}
\usepackage{amsthm}
\usepackage{amscd}
\usepackage{amsfonts}
\usepackage{ascmac}
\usepackage[usenames]{color}
\usepackage{graphics}
\usepackage{DotArrow}

\theoremstyle{plain}

\newtheorem{theorem}{Theorem}
\newtheorem{proposition}[theorem]{Proposition}

\newtheorem{definition}[theorem]{Definition}

\newtheorem{remark}[theorem]{Remark}

\newcommand{\bfC}{{\mathbf C}}

\newcommand{\bfP}{{\mathbf P}}

\newcommand{\barj}{{\overline j}}

\newcommand{\barl}{{\overline \ell}}

\newcommand{\barpsi}{{\overline \psi}}

\newcommand{\barpartial}{{\overline \partial}}

\newcommand{\mapright}[1]{\smash{\mathop{   \hbox to 0.7cm{\rightarrowfill}}
  \limits^{#1}}}

\newcommand{\Ker}{{\rm Ker}}

\newcommand{\Ka}{K\"ahler }

\newcommand{\Fut}{\mathrm{Fut}}

\newcommand{\grad}{\mathrm{grad}}

\newcommand{\Aut}{\operatorname{Aut}}
\newcommand{\Isom}{\operatorname{Isom}}

\title[Weighted extremal K\"ahler metrics]{On weighted extremal K\"ahler metrics}
\author[]{Akito Futaki}
\address{Yau Mathematical Sciences Center, Tsinghua University, Haidian District, 
Beijing 100084, China}
\email{futaki@tsinghua.edu.cn}\thanks{}

\begin{document}

\begin{abstract}
The notion of weighted extremal K\"ahler metrics extends the classical notion of Calabi's extremal \Ka metrics, but includes many well-studied objects in K\"ahler geometry such as K\"ahler-Ricci solitons and Sasaki-Einstein metrics. In this paper, after explaining how this notion grew out, we will try to survey recent works concerning the YTD conjecture on weighted extremal K\"ahler metrics.
\end{abstract}

\maketitle


\section{Introduction}
This manuscript is based on the author's talk at ICCM 2025 held in Shanghai in January 2026.

In K\"ahler geometry, after the celebrated works of Yau solving the Calabi conjecture (\cite{yau77}, \cite{yau78}, also \cite{aubin76}) there have been extensive studies to prove 
the equivalence between the
existence of special K\"ahler metrics and the possession of an algebraic stability property called the K-(poly)stability 
known as the (generalized) Yau-Tian-Donaldson conjecture (\cite{yau87}, \cite{tian97}, \cite{donaldson02}).
The special metrics are typically,
\begin{itemize}
\item[-] K\"ahler-Einstein (KE for short) metrics, 
\item[-] constant scalar curvature K\"ahler (cscK for short) metrics including KE metrics,
\item[-] extremal K\"ahler metrics including cscK metrics, 
\item[-] K\"ahler-Ricci solitons,
\item[-] Mabuchi solitons,
\item[-] Sasaki-Einstein metrics, 
\item[-] conformally K\"ahler, Einstein-Maxwell (cKEM for short)metrics, 
\item[-] $\mu$-scalar curvature considered in \cite{Inoue19}, \cite{Inoue22} including cKEM metrics.
\end{itemize}


In this paper we discuss on special K\"ahler metrics, called weighted extremal K\"ahler metrics, which include all the above examples.
This notion was introduced by Lahdili \cite{Lahdili18}
extending the notion of conformally K\"ahler, Einstein-Maxwell equations. 
The new notion allows us a unified treatment of those various types canonical K\"ahler metrics.

In section 2, we will review the notion of weighted extremal metrics.

In section 3, we will review conformally K\"ahler, Einstein-Maxwell equations. 

In section 4, we show many standard results on extremal K\"ahler metrics going back to Calabi \cite{calabi85} extend to 
weighted extremal K\"ahler metrics.

In section 5, we will review the K-polystability characterization of weighted solitons. We mainly summarize the works 
of Berman-Boucksom-Jonsson \cite{BBJ} on uniform D-stability and Li \cite{LiChi22} on $G$-uniform K or D-stablity
for Fano K\"ahler-Einstein manifolds.

In section 6, we will review the works of Li \cite{LiChi22b}, Boucksom-Jonsson \cite{BoucJon25} and Darvas-Zhang \cite{DarvasZhang25} on the cscK and weighted extremal  metrics.

In section 7, we review some works on conformally K\"ahler, Einstein-Maxwell equations. 

\section{Weighted extremal metrics}

Weighted scalar curvature, or $(v,w)$-scalar curvature, was defined by Lahdili \cite{Lahdili18}.
Let $X$ be a compact connected K\"ahler manifold of dimension $n$ and $\Omega$ its K\"ahler class.
Denote by $\Aut(X)$ the automorphism group of $X$, and by $\mathfrak h(X)$
 the Lie algebra of $\Aut(X)$. Then $\mathfrak h(X)$ consists of holomorphic vector fields on $X$.
Let $\mathfrak h_{red} \subset \mathfrak h(X)$ be the subalgebra of $\mathfrak h(X)$
consisting of holomorphic vector fields $\xi$ which are written as 
$$\xi=\grad'u= \sum_{i,j}g^{i\bar{j}}\frac{\partial u}{\partial z^{\bar{j}}}\frac{\partial}{\partial z^i}$$
for some complex valued smooth function $u \in C^\infty_{\bfC}(X)$ on $X$ and for any K\"ahler form 
$\omega = \sum_{i,j} \sqrt{-1} g_{i\bar{j}} dz^i \wedge dz^{\bar{j}} \in \Omega$. $\mathfrak h_{red}$ is called 
the reduced subalgebra. We will denote its Lie subgroup in $\Aut(X)$ by $\Aut_{red}(X)$. Note that 
the last condition is equivalent to $i(\xi) \omega = \sqrt{-1}\,\barpartial u$.
When $\Omega = c_1(L)$ for some ample line bundle $L \to X$, $\mathfrak h_{red}$ is
the Lie algebra of $\Aut(X,L)$. 


Let $T$ be a compact real torus in $\Aut(X)$ such that $\mathfrak t \otimes \mathbf C \subset \mathfrak h_{red}$,
and $\omega \in \Omega$ a $T$-invariant K\"ahler form. We often identify $\omega$ with its K\"ahler metric $g=(g_{i\barj})$.
Then $T$ acts on $(X,\omega)$ in the Hamiltonian way. 
Let $\mu_\omega : X \to \mathfrak t^\ast$ be the moment map. 
Then $\Delta:=\mu_\omega(X)$ is a compact convex polytope.
The moment polytope $\Delta$ is independent of $\omega \in \Omega$ up to translation.


Let $v$ be a positive smooth function on $\Delta$. 
Write $v = v(\mu)$ as a function on $\Delta$ and
$v(\mu_\omega) := \mu_\omega^\ast v$ as a positive smooth function on $X$.

Define $v$-scalar curvature $S_v(\omega)$ of a $T$-invariant K\"ahler form $\omega$ by
$$ S_v(\omega):= v(\mu_\omega) S(\omega) + 2 \Delta_\omega v(\mu_\omega) 
+ \langle g_\omega,\mu_\omega^\ast Hess(v)\rangle$$
where $S(\omega)$ denotes the 
scalar curvature 
$$ S(\omega) = - g^{i\barj} \frac{\partial^2 }{\partial z^i \partial z^\barj} \log \det (g_{l\barl})$$
of $\omega$ and
$$\langle g_\omega,\mu_\omega^\ast Hess(v)\rangle = g^{i\barj} v_{\alpha\beta}\mu^\alpha_i\mu^\beta_\barj$$
where $Hess(v)$ of $v$ is the Hessian of $v$ on $\frak t^\ast$ and
$\mu_\omega(p) = (\mu^1(p), \cdots, \mu^\ell(p))$ with $d\mu^\alpha = i(\xi^\alpha)\omega$ for a basis $\xi^1, \cdots, \xi^\ell$ of $\frak t$.


Let $w$ be another positive smooth function on $\Delta$.
Define $(v,w)$-scalar curvature $S_{v,w}$ by
$$ S_{v,w} = \frac{S_v}{w(\mu_\omega)}.$$
We say $g$ is a weighted constant scalar curvature K\"ahler (cscK for short) metric if $S_{v,w}$ is constant.
As we explain in the next section, the notion of $(v,w)$-scalar curvature was originally introduced as a generalization of
conformally K\"ahler, Einstein-Maxwell metrics.
Later it turned out that 
the $(v,w)$-cscK metrics include much more unexpected examples as above.


We call $g$ a weighted extremal metric or $(v,w)$-extremal metric if 
$$\grad'S_{v,w} = g^{i\barj}\frac{\partial S_{v,w}}{\partial z^\barj}\frac{\partial}{\partial z^i}$$
is a holomorphic vector field. 
Then $g$ satisfies
$$ S_{v,w} = \ell\circ \mu_\omega$$
for some affine function $\ell$ on $\Delta$. Putting $\tilde{w} = w\ell$, we could write $S_v = \tilde{w}$. 
Though $\tilde{w}$ may no longer be positive, the equation is the same form of $S_{v,\tilde{w}}$ being constant.

A slightly different formulation of weighted extremal metrics is also given in \cite{BJT26}. The result of Chen-Cheng \cite{CC21b}
that the properness of Mabuchi energy is equivalent to the existence of cscK metrics was extended to
the weighted case by Han-Liu \cite{HL25a} and Di Nezza-Jubert-Lahdili \cite{DJL25}.

\section{Conformally K\"ahler Einstein-Maxwell metrics}
Lahdili's introduction of weighted scalar curvature was motivated by
conformally, K\"ahler, Einstein-Maxwell metrics previously considered by LeBrun \cite{LeBrun16} and Apostolov-Maschler \cite{AM}.

Let $(X,J)$ be a compact K\"ahler manifold where $J$ is the integrable almost complex structure.
A Hermitian metric (i.e. $J$-invariant Riemannian metric) 
$\tilde{g}$ on $(X,J)$ is called 
a conformally K\"ahler, Einstein-Maxwell metric
(cKEM metric for short) if it satisfies the following three conditions:

(a) There exists a positive smooth function $f$ on $X$ 
such that 
$g=f^2\Tilde{g}$  is K\"ahler.

(b) The Hamiltonian vector field $K=J\mathrm{grad}_gf$
 is Killing for both $g$ and $\Tilde{g}$.

(c)  The scalar curvature $s_{\Tilde{g}}$ of $\Tilde g$ is constant.


\noindent
The condition (c) is equivalent to 
\begin{equation*}\label{eq:3.1}
s_{\Tilde{g}}=2\left(\dfrac{2m-1}{m-1}\right)f^{m+1}\Delta_g
\left(\frac{1}{f}\right)^{m-1}+s_gf^2=\text{const},
\end{equation*}
where $\dim_{\mathbf{C}}X=m$.
Putting $f = (\langle \xi, \mu_\omega \rangle +a )$,
$v(p) = (\langle \xi, p \rangle + a)^{-2m+1}=f^{-2m+1}$,
$w(p) = \langle \xi, p \rangle +a)^{-2m-1}=f^{-2m-1}$ 
for $\xi \in \mathfrak t$ and $a \in \mathbf R$ large enough positive constant,
one can compute
$$ s_{\Tilde{g}} = S_{v,w}.$$

In view of K\"ahler geometry, rather than trying to find a Hermitian metric $\tilde g$,  it will be more natural
to try to find a K\"ahler metric $g$ and a Killing potential $f$ such that the conformal Hermitian metric 
$\tilde g = f^{-2} g$ has constant scalar curvature. In this regard, we may call $g$ a conformally Einstein-Maxwell
K\"ahler metric (\cite{FO_reductive17}).

%
If $\Tilde g$ is a cKEM metric and if $\dim_{\mathbf R} = 4$, then one obtains a solution 
$(X, h, F)$ of the following Einstein-Maxwell equation studied in General Relativity 
(see LeBrun \cite{L1}):

(i) $h$ is a Riemannian metric. (In our case $h = \Tilde g$).

(ii) $F$ is a real 2-form.

(iii) $dF = 0$, $d\ast F=0$, $[\mathrm{Ric} + F\circ F]_0 = 0$. Here $(F\circ F)_{jk} = 
F_j{ }^\ell F_{\ell k}$.
(In our case $F^+$ is the K\"ahler form 
$\omega_g$, and $F^- = \frac12 f^{-2} \rho_0(\Tilde g)$ with $\rho_0(\Tilde g)$ the traceless Ricci form of 
$\Tilde g$.)
%
\section{Standard results on extremal K\"ahler metrics and weighted extremal K\"ahler metrics}
Note that when $v=w=1$ the weighted scalar curvature $S_{v,w}$ is the usual scalar curvature $S$ of $g$.
Many standard results for cscK problem extend to the weighted cscK problem as follows.
Define $L_v\varphi$ for complex valued smooth functions $\varphi$ by
$$ L_v \varphi = \nabla^i\nabla^j(v(\mu_\omega)\nabla_i\nabla_j \varphi),$$
and call $L_v$ the $v$-twisted Lichnerowicz operator. Obviously, 
$$ \int_X (L_v\varphi)\, \barpsi\, \omega^m =  \int_X \varphi\, \overline{L_v \psi}\, \omega^m,$$
and 
$\Ker\, L_v$ is isomorphic to $\frak h_{red}$. 
We also define $L_{v,w}$ by
$$ L_{v,w} = \frac 1{w(\mu_\omega)}\, L_v.$$
Consider $g_{ti\barj} = g_{i\barj} + t\varphi_{i\barj}$. Then
\begin{equation*}\label{Lich1}
\left.\frac d{dt}\right|_{t=0} S_v(g_t) = - L_v \varphi + S_v^i\,\varphi_i,
\end{equation*}
\begin{equation*}\label{Lich2}
\left.\frac d{dt}\right|_{t=0} S_{v,w}(g_t) = - L_{v,w} \varphi + S_{v,w}^i\,\varphi_i.
\end{equation*}
Using this, by the same argument as in the usual unweighted scalar curvature case, we can show the following.
\begin{proposition}\label{Lich3}
A critical point of the weighted Calabi functional 
$$g \mapsto \int_X S_{v,w}^2 (g)\,w(\mu_{\omega})\, \omega^m$$ 
is a weighted extremal metric.
\end{proposition}


\begin{proposition}\label{Lich3.1}
Let $h_\xi \in \Ker\,L_v$ be the real Killing potential of $\xi \in \frak t$, 
i.e. $\sqrt{-1}\,\grad'h_\xi = \xi'$.
Then $\Fut_v$ and $\Fut_{v,w}$ defined by
\begin{equation}\label{Lich4}
 \Fut_v (\xi) = \int_X (S_v - c_v)\, h_\xi\,\omega^m
 \end{equation}
 and 
 \begin{equation}\label{Lich5}
 \Fut_{v,w} (\xi) = \int_X (S_{v,w} - c_{v,w})\, h_\xi\,w(\mu_\omega)\,\omega^m
 \end{equation}
are independent of choice of $\omega \in \Omega$ where 
$c_{v,w} = \int_X S_v\, \omega^m/\int_X w(\mu_\omega)\, \omega^m$ 
and $c_{v} = c_{v,1}$
which
are independent of $\omega \in \Omega$. 
In particular, 
if $ \Fut_{v,w} \ne 0$ then $X$ does not admit a $(v,w)$-cscK metric in $\Omega$.
\end{proposition}


\begin{remark}\label{Lich5.1}If $g$ is a $(v,w)$-extremal metric with non-constant $S_{v,w}$ then
$$ \Fut_{v,w}(J\grad S_{v,w}) = \int_X (S_{v,w} - c_{v,w})^2\, w(\mu_\omega)\,\omega^m > 0.$$
\end{remark}
\begin{remark}Calabi's decomposition theorem for extremal K\"ahler metrics holds for weighted extremal K\"ahler metrics.
In particular, as an extension of the theorem of Lichnerowicz-Matsushima, it follows that the existence of weighted cscK implies 
$$\Aut_{red}^T(X)_{id} = (\Isom_{red}^T(X)_{id})^\bfC$$
where $\Aut_{red}^T(X)_{id} $ is the identity component of the centralizer of $T$ in $\Aut_{red}(X)$ and $\Isom_{red}^T(X)_{id}$ is the subgroup of $\Aut_{red}^T(X)_{id}$ consisting of isometries of $g$ in $\Aut_{red}^T(X)_{id}$.
This shows that $\Aut_{red}^T(X)_{id}$ is reductive.
\end{remark}

\section{K-polystability characterization for KE metrics and weighted solitons}
As for $K$-polystability characterization of K\"ahler-Einstein metrics, 
Chen-Donaldson-Sun \cite{CDS3} and Tian \cite{Tian12} proved that a Fano manifold $X$ admits a K\"ahler-Einstein metric if 
$\left(X, K_X^{-1}\right)$ is K-polystable (with only if part due to \cite{donaldson05}, \cite{stoppa0803}, \cite{BDL}, \cite{Berman16}). Datar-Szekelyhidi \cite{DatarSzeke16} and Chen-Sun-Wang \cite{CSW} also gave other proofs. 
All these proofs use the Cheeger-Coclding-Tian theory on Gromov-Hausdorff convergence under Ricci curvature bound. The definitions of the basic keywords are as follows.

\begin{definition}\label{tc} Let $(X,L)$ be a polarized complex manifold, i.e. a smooth projective complex variety $X$ with an ample line bundle $L$. An ample test configuration $(\mathcal{X},\mathcal{L})$ for $(X,L)$ consists of the following data:
\begin{enumerate}
\item[(i)] a flat, $\mathbb{C}^*$-equivariant projective morphism of normal schemes $\pi: \mathcal{X} \to \mathbb{P}^1$;
\item[(ii)] a $\mathbb{C}^*$-equivariant relatively ample $\mathbb{Q}$-line bundle $\mathcal{L}$;\\
\item[(iii)] $(\mathcal{X}_t,\mathcal{L}_t) \simeq (X,L)$ for $t\neq 0$ where $\mathcal X_t = \pi^{-1}(t)$ and $\mathcal L_t = \mathcal L\vert_{\mathcal X_t}$.
	\end{enumerate}
A test configuration is called product (resp. product type) if it is of the form $(\mathcal X, \mathcal L)$ (resp. $(\mathcal X, \mathcal L + c\mathcal X_0)$) where 
$(\mathcal X,\mathcal L)\vert_{\mathbf P^1 \backslash \{\infty\}}$ is $\mathbf C^\ast$-equivariantly isomorphic to $(X,L)\times \mathbf C$ 
for some one parameter subgroup $\mathbf C^\ast \subset \Aut(X,L)$.
A test configuration is called trivial if it is of the form $(\mathcal X, \mathcal L + c\mathcal X_0)$ where 
$(\mathcal X,\mathcal L)$ is $\mathbf C^\ast$-equivariantly isomorphic to $(X,L)\times \mathbf \mathbf P^1$ 
for some one parameter subgroup $\mathbf C^\ast \subset \Aut(X,L)$.
\end{definition}

\begin{definition}\label{df} Donaldson-Futaki invariant $DF(\mathcal X, \mathcal L)$ of the test configuration $(\mathcal X, \mathcal L)$ is defined by
$$ DF(\mathcal X, \mathcal L) = \frac{(K_{\mathcal X/\mathbb P^1}\cdot \mathcal L^n)}{V} + \bar S \frac{(\mathcal L^{n+1})}{(n+1)V}$$
where $V = L^n$ the volume, $\bar S = - n V^{-1} (K_X\cdot L^n)$ the average scalar curvature, and $K_{\mathcal X/\mathbf P^1}= K_{\mathcal X} - \pi^\ast (K_{\mathbf P^1})$ the relative canonical line bundle. Here we use the intersection formula of \cite{Odaka13}, \cite{wangxw12}, \cite{LX14}.
\end{definition}

\begin{definition} 
\begin{enumerate}
\item[(i)] A polarized manifold $(X,L)$ is said to be K-semistable if $DF(\mathcal X, \mathcal L) \ge 0$ for all ample test configurations.
\item[(ii)] 
A polarized manifold $(X,L)$ is said to be K-polystable if it is K-semistable and the equality $DF(\mathcal X, \mathcal L) = 0$
holds if and only if $(\mathcal X, \mathcal L)$ is product type.

\item[(iii)] A polarized manifold $(X,L)$ is said to be stable if $DF(\mathcal X, \mathcal L) > 0$ for all non-trivial ample test configurations.
\end{enumerate}
\end{definition}

For Fano manifolds, these definitions are applied for $L = - K_X$.

On the other hand, 
Berman-Boucksom-Jonsson \cite{BBJ} proved by a variational method without using the Gromov-Hausdorff convergence that a Fano manifold $X$ admits a unique K\"ahler-Einstein metric if and only if
$X$ satisfies uniform $K$-stability. 
In this case the automorphism group is discrete.
Here, the definitions of basic key words are as follows. 
Let 
$\omega$ be a fixed K\"ahler form, and 
$$
\mathcal{H}_\omega=\left\{v \in C^{\infty}(X) \mid \omega_v:=\omega+i\partial \bar{\partial} v>0\right\}
$$
the space of K\"ahler forms cohomologous to $\omega$.

\begin{definition}
The Ding functional is a map $D : \mathcal H \to \mathbb R$ defined by
$$D(u) := L(u) - E(u)$$
where
$$ L(u) := -\log \int_X \exp(-u + f) \omega^n$$
with $Ric(\omega) = \omega + i\partial\barpartial f$, 
and
$$
E(u)=\frac{1}{(n+1) V} \int_X u\, \sum_{j=0}^n\,\omega_u^j \wedge \omega^{n-j} .
$$
The functional $E$ is called the Monge-Amp\'ere energy (a.k.a. Aubin-Mabuchi energy). 
\end{definition}
\begin{definition}
The Ding functional is said to be coercive if there exists an $\epsilon > 0$ such that
\begin{equation}\label{coer}
 D(u) \ge \epsilon J(u) - \frac1{\epsilon}
\end{equation}
where
$$ J(u) := \frac1V \int_X\, u\, \omega^n - E(u).$$
The functional $J \ge 0$ is called Aubin's $J$-functional, and plays the role of the norm on $\mathcal H$. 
\end{definition}
$D$ is convex along plurisubharmonic (psh) geodesics \cite{Berndt09}. If $u \in \mathcal H$ is an absolute minimum 
of $D$ then $\omega_u$ is a K\"ahler-Einstein metric. On $\mathbf P^1(\mathbf C)$, the inequality \eqref{coer} is the Sobolev inequality of exceptional type. 

For an ample test configuration, we choose a metric of $\mathcal L$ which gives a ray $h_s = e^{-\phi_s}$ of metrics of $L_t \to X_t$
over $t \in \mathbf P^1$ where $s = - \log |t|$. The coerciveness \eqref{coer} implies the inequality of the slopes at $s = \infty$ :
\begin{equation} 
\lim_{s \to \infty} \frac{D(\phi_s)}{s} \ge \epsilon \lim_{s\to \infty}\frac{J(\phi_s)}{s}.
\end{equation}
But one can show the slopes indeed exist and can be computed as follows. 
$$ \lim_{s\to \infty}\frac1s E(\phi_s) = \frac{\mathcal L^{n+1}}{(n+1)V} =: E_{na}(\mathcal X,\mathcal L),$$
$$ \lim_{s\to \infty}\frac1s L(\phi_s) = \mathrm{lct}(\mathcal X,\mathcal D_{\mathcal X,\mathcal L};\mathcal X_0) -1 =: 
L_{na}(\mathcal X,\mathcal L)
$$
where $\mathcal D_{\mathcal X,\mathcal L} \sim_\mathbb Q - \mathcal L - K_{\mathcal X/\mathbb P^1}$ supported on $\mathcal X_0$, and 
$$ \mathrm{lct}(\mathcal X,\mathcal D_{\mathcal X,\mathcal L};\mathcal X_0) = \sup \{t\ |\ (\mathcal X, \mathcal D_{\mathcal X,\mathcal L} + t \mathcal X_0) \text{\ sublog canonical}\} $$
is the log-canonical threshold. 
Here sublog canonical means that $$\mathrm{coeff}_E (K_\mathcal N - \mu^\ast (K_\mathcal X + \mathcal D_{\mathcal X,\mathcal L} + t \mathcal X_0) > -1$$
for  any prime divisor $E$ over $\mathcal X$ (i.e. $\mu : \mathcal N \to \mathcal X$ is proper birational, $\mathcal N$ regular and $E \subset \mathcal N$).
Hence
\begin{eqnarray*}
 \lim_{s\to \infty}\frac1s D(\phi_s)  &=&  \mathrm{lct}(\mathcal X,\mathcal D_{\mathcal X,\mathcal L};\mathcal X_0)  - 1 
- \frac{\mathcal L^{n+1}}{(n+1)(-K_X)^n}  \\
&=&L_{na}(\mathcal X,\mathcal L) - E_{na}(\mathcal X,\mathcal L) =: D_{na}(\mathcal X,\mathcal L).
\end{eqnarray*}
For $J$-functional, replacing $\mathcal X$ by a birational model dominating both $\mathcal X$ and $X\times \mathbf P^1$, we define
 $$ J_{na}(\mathcal X, \mathcal L) := \frac{\mathcal L\cdot \rho^\ast (L\times \bf P^1)^n }{L^n} - \sum \frac{\mathcal L^{n+1}}{(n+1)L^n}
$$
with $\rho : \mathcal X \to X \times \mathbb P^1$. Then
$$ \lim_{s\to \infty}\frac1s J(\phi_s) = J_{na}(\mathcal X, \mathcal L).$$

\begin{definition}
A Fano manifold $X$ is said to be uniformly D-stable if there is an $\epsilon > 0$ such that
$$ D_{na}(\mathcal X, \mathcal L) \ge \epsilon J_{na}(\mathcal X, \mathcal L)$$
for any test configuration $(\mathcal X, \mathcal L)$ of $(X, -K_X)$.
\end{definition}
Uniform stability describes the behavior of slopes of Ding-energy and $J$-functional, which ensures that $D$ is proper, i.e. $D \to \infty$ when $J \to \infty$. Thus, the uniform stability is an asymptotic description of coercivity or properness of the Ding functional.

The result of Berman-Boucksom-Jonsson \cite{BBJ} is that, on a Fano manifold $X$, there exists a unique K\"ahler-Einstein metric 
if and only if $X$ is uniformly Ding-stable. Further they extended this result for the problem of finding twisted K\"ahler-Einstein metrics. It is known that uniform Ding-stability is equivalent to uniform K-stability where the latter use the non-Archimedean Mabuchi functional $M_{na}$, see below.
 \begin{definition}
 For $u \in \mathcal{H}_\omega$ we define the Mabuchi $K$-energy functional $M(u)$ by
$$
	M(u):=\frac{1}{V} \int_0^1 d t \int_X \dot{v}_t\left(S\left(\omega_{v_t}\right)-\underline{S}\right) w_{v_t}^n
	$$
where $V=\int_X \omega^n$, $\underline{S}=\frac{1}{V} \int_X S(\omega)\, \omega^n$ is the average of the scalar curvature, and 
	$\left\{v_{t} \mid 0 \leq t \leq 1\right\}$ is a smooth path between $0$ and $u$
	such that $\omega_{v_{t}}>0$.
	Mabuchi \cite{mabuchi87} showed that $M(u)$ is independent of the path $v_t$.
\end{definition}
Mabuchi energy has a decomposition, called Chen-Tian formula, 
	$$
		M(u)=H(u)+\underline{S} E(u)-n E_{\operatorname{Ric}} (u), $$
where 
$H$, called the entropy, is defined by
$$	
		H(u)=V^{-1} \int_X \log \left(\omega_u^n / \omega^n\right) \omega_u^n,$$
$E$, called the Monge-Ampere energy is defined by
$$
\operatorname{E}_\omega (u)=\operatorname{E}(u)=\frac{1}{(n+1) V} \sum_{j=0}^n \int_X u\, \omega_u^j \wedge \omega^{n-j} 
$$
and $E_{\operatorname{Ric}(\omega)}$ is defined by 
$$
		E_{\operatorname{Ric}}(u)=\frac{1}{n V} \sum_{j=1}^{n} \int_X u\, \operatorname{Ric}(\omega) \wedge \omega_u^{j-1} \wedge \omega^{n-j}.
		$$
See section 3 in \cite{Chen00} for its derivation.
\begin{definition}
The Mabuchi energy $M$ is said to be proper if $M(u) \to \infty$ when $J(u) \to \infty$. 
$M$ is said to be coercive if there exists an $\epsilon > 0$ such that
\begin{equation}\label{Mcoer}
 M(u) \ge \epsilon J(u) - \frac1{\epsilon}.
 \end{equation}
\end{definition}

As in the case of Ding functional, 
for an ample test configuration, we choose a metric of $\mathcal L$ which gives a ray $h_s = e^{-\phi_s}$ of metrics of $L_t \to X_t$
over $t \in \mathbf P^1$ where $s = - \log |t|$. The coerciveness \eqref{Mcoer} implies the inequality of the slopes at $s = \infty$ :
\begin{equation} 
\lim_{s \to \infty} \frac{M(\phi_s)}{s} \ge \epsilon \lim_{s\to \infty}\frac{J(\phi_s)}{s}.
\end{equation}
Then 
\begin{equation} 
\lim_{s \to \infty} \frac{M(\phi_s)}{s} = M_{na}(\mathcal X, \mathcal L) 
\end{equation}
where $M_{na}$ is computed as follows (\cite{BHJ19}, \cite{Dyre18}).

\begin{definition}The non-Archimedean Mabuchi functional is defined by
\begin{equation}\label{mna} 
M_{na}(\mathcal X, \mathcal L)  = V^{-1} (K^{\log}_{\mathrm X/\mathbf P^1}\cdot \mathcal L^n) + \frac{\bar S}{(n+1)V} \mathcal L^{n+1}
\end{equation}
where
$$ K^{\log}_{\mathcal X/\mathbf P^1} := K^{\log}_{\mathcal X} - \pi^\ast K^{\log}_{\mathbf P^1} := K_{\mathcal X} + {\mathcal X}_{0,\mathrm{red}} - \pi^\ast (K_{\mathbf P^1} + \{0\}),
$$
which is called the relative logarithmic canonical line bundle. 
\end{definition}
Note that $DF$ and $M_{na}$ are related by
$$ DF(\mathcal X, \mathcal L) = M_{na} (\mathcal X, \mathcal L) + V^{-1} (\mathcal X_0 - \mathcal X_{0,\mathrm{red}})\cdot \mathcal L^n \ge M_{na} (\mathcal X, \mathcal L).$$

As remarked above, in this case of Berman-Boucksom-Jonsson \cite{BBJ}, the automorphism group of $X$ is discrete. Li \cite{LiChi22} extended their result for non-discrete automorphism case. 
\begin{definition}\label{Gtest}
Let $G$ be a reductive subgroup of the
identity component $\Aut(X,L)_0$ of $\Aut(X,L)$.
A test configuration $(\mathcal X,\mathcal L)$ of a polarized manifold $(X,L)$ is a $G$-equivariant test configuration
if $G$ acts on $(\mathcal X,\mathcal L)$ commuting with the $\mathbf C^\ast$-action of $(\mathcal X,\mathcal L)$ and the $G$-action on 
$(\mathcal X,\mathcal L)\times_{\mathbf C} \mathbf C^\ast \cong (X,L)\times \mathbf C^\ast $ coincides with the action
of $G$ on the first factor of $(X,L)\times \mathbf C^\ast $. 
\end{definition}
\begin{definition} Let $T$ be the identity component of the center of $G$.
The reduced $J$-norm $J_{T,na}(\mathcal X,\mathcal L)$ is defined by
$$ J_{T,na}(\mathcal X,\mathcal L) = \inf_{\xi \in N_{\mathbf R}} J_{na}(\mathcal X_\xi,\mathcal L_\xi)$$
where $N_{\mathbf Z} = \mathrm{Hom}(\mathbf C^\ast,T)$, $N_{\mathbf R} = N_{\mathbf Z}\otimes \mathbf R$
and $(\mathcal X_\xi,\mathcal L_\xi)$ is the $\xi$-twist of of $(\mathcal X,\mathcal L)$, i.e. the test configuration whose $\mathbf C^\ast$-action is twisted
by $exp(t\xi)$ against the $\mathbf C^\ast$-action of $(\mathcal X,\mathcal L)$. 
\end{definition}

\begin{definition}\label{G-uni}A Fano manifold is $G$-uniformly $K$-stable if 
if there is an $\epsilon > 0$ such that
$$ M_{na}(\mathcal X, \mathcal L) \ge \epsilon J_{T,na}(\mathcal X, \mathcal L)$$
for any test configuration $(\mathcal X, \mathcal L)$ of $(X, -K_X)$.
Replacing $M_{na}$ by $D_{na}$, one defines $G$-uniform D-stability.
\end{definition}

Li \cite{LiChi22} proved that \\
(i)\ \ a Fano manifold $X$ admits a K\"ahler-Einstein metric if $X$ is $G$-uniformly K-stable, 
 or equivalently, $G$-uniformly D-stable, and that \\
 (ii)\ \ if $G$ contains the maximal torus of $\Aut(X,L)$ then
 the converse also holds. 
 
On the other hand, Liu-Xu-Zhuang \cite{LXZ22} proved that\\ 
(iii)\ \ on Fano manifolds, $G$-uniform $K$-stability is equivalent to K-polystability if $G$ contains the maximal torus. 

Hence, these works give another proof of K-polystability characteriztion of the existence of K\"ahler-Einstein metrics on Fano manifolds.

Further,  Han-Li \cite{HanLi23}, 
Li \cite{LiChi21} and Blum-Liu-Xu-Zhuang \cite{BLXZ} extended these results to obtain K-stability characterization for
weighted solitons on Fano manifolds. Here weighted solitons are special cases of weighted extremal metrics
which we discussed above. Typical weighted solitons are K\"ahler-Ricci solitons. 

\section{Modified K-stability characterizations for cscK metrics and weighted extremal metrics}
Now we turn to the case of cscK and weighted extremal metric case. As an analytic characterization, Chen-Cheng \cite{CC21b} proved that for a polarized manifold
$(X,L)$, the properness of the Mabuchi energy
is equivalent to the existence of constant scalar curvature metric. 

As an algebraic characterization, Li \cite{LiChi22b} proved that on a polarized manifolds $(X,L)$,
a cscK metric exists if $(X,L)$ is $G$-uniformly K-stable for models. Below are the basic key definitions.

\begin{definition}
For a polarized manifold $(X,L)$, \\
(i)\ \ a pair $(\mathcal X,\mathcal L)$ is called a model of $(X,L)$ if $(\mathcal X,\mathcal L)$
satisfies all the conditions in the Definition \ref{tc} of test configurations except the ampleness of $\mathcal L$,\\
(ii)\ \ $(X,L)$ is said to be $G$-uniformly K-stable for models if $(X,L)$ satisfies all the inequality in Definition \ref{G-uni}
are satisfied for all models of $(X,L)$.
\end{definition}

Since models are larger class than test configurations, the condition (ii) above is stronger than the condition in Definition
\ref{G-uni}. Recall in Definition \ref{Gtest}, $G$ is a reductive subgroup of the
identity component $\Aut(X,L)_0$ of $\Aut(X,L)$. 
Boucksom-Jonsson \cite{BoucJon25} proved that the converse of Li's result holds when $G = \Aut_0(X,L)$.
\begin{theorem}[\cite{LiChi22b}, \cite{BoucJon25}]
Let $(X, L)$ be a polarized K\"ahler manifold. 
The following conditions are equivalent:\\
(i) there exists a cscK metric $\omega \in c_1(L)$;\\
(ii) $G := \Aut_0(X,L)$ is reductive and $(X,L)$ is $G$-uniformly K-stable with respect to models.
\end{theorem}
There is no exact statement of the above result in the paper \cite{BoucJon25}, but it is a combination of Theorem B and Lemma 8.17 in \cite{BoucJon25}. See also Theorem 1.3 in \cite{DarvasZhang25}.

Boucksom-Jonsson \cite{BoucJon25} describe their result in terms of non-Archimedean pluripotential theory.
In the (Archimedean) pluripotential theory by Bedford-Taylor \cite{BedfordTaylor},  
for non-smooth plurisubharmonic (PSH for short) functions $v$, the wedge product of currents
$i\partial\overline{\partial}v \wedge \cdots \wedge i\partial\overline{\partial}v$ makes sense. 
Similarly, we may consider $\omega$-plurisubhamonic functions $u$, and put
$$
\begin{aligned}
	& \mathcal E=\left\{u \in \operatorname{PSH}(X, \omega) \mid \int_X \omega_u^n=V\right\} \\
	& \mathcal E^{1}=\left\{u \in \mathcal E \mid \int_X |u|\, \omega_u^n\leq \infty\right\}.
\end{aligned}
$$
The spaces $\mathcal E$ and $\mathcal E^{1}$ were introduced by Guedj-Zeriahi \cite{GuedjZeriahi}. 
Boucksom-Eyssidieux-Guedj-Zeriahi \cite{BEGZ10} and Berman-Boucksom-Guedj-Zeriahi \cite{BBGZ13} extended
the functional $E$ to $\mathcal E^1$. Darvas \cite{Dar17} defined $d_1$-distance and proved
that the completion of $\mathcal H_{\omega}$ with respect to $d_1$ is $\mathcal E_1$ and that the functional $E$ is $d_1$-continuous.
Boucksom-Jonsson \cite{BoucJon25} consider $\mathcal H_{na}$ to be the space of 
all ample test configurations of $(X,L)$ 
equipped with a naturally defined Darvas metric $d_{1,na}$, and $\mathcal E^1_{na}$ to be the metric completion of $\mathcal H_{na}$
with respect to  $d_{1,na}$. $M_{na}$ extends to $\mathcal E^1_{na}$ as a lower semi-continuous functional. They define a polarized manifold $(X,L)$ to be\\
(a)\ \ $\hat K$-polystable if $M_{na}(\varphi)\ge 0$ for all $\varphi \in \mathcal E^1_{na}$,
with the equality only if $\varphi$ is a real product test configuration.\\
(b)\ \ uniformly $\hat K$-polystable if $\Aut_0(X,L)$ is reductive and there exists $\epsilon > 0$ such that
$M_{na}(\varphi) \ge \epsilon\,d_{1,na}(\varphi, \mathcal P_{\mathbf R})$ for all $\varphi \in \mathcal E^1_{na}$ 
where $\mathcal P_{\mathbf R}$ is the subspace of real product configurations.
For the cscK problem, they prove the following.
\begin{theorem}[Boucksom-Jonsson \cite{BoucJon25}] For a polarized manifold $(X,L)$, the following are equivalent.
\begin{enumerate}
\item[(i)]  There exists a cscK metric in $c_1(L)$;
\item[(ii)] $(X,L)$ is $\hat K$-polystable;
\item[(iii)] $(X,L)$ is uniformly $\hat K$-polystable.
\end{enumerate}
\end{theorem}
This result is extended by Boucksom-Jonsson \cite{BoucJon25} for weighted extremal metrics
as follows. 
Assume that $v$ is log-concave. Then the following are equivalent:\\
(i) there exists a $T$-invariant weighted extremal metric $\omega \in c_1(L)$;\\
(ii) $(X,L)$ is $T$-equivariantly relatively weighted $\hat K$-polystable;\\
(iii) $(X,L)$ is uniformly $T$-equivariantly relatively weighted $\hat K$-polystable.\\
We do not describe the precise definitions of the terminologies in the above result, and
just refer to Theorem B in their paper. 
See also related works \cite{HL25b}, \cite{Hashimoto26}. 

Darvas-Zhang \cite{DarvasZhang25} also gave a necessary and sufficient condition for the existence
of cscK metrics by an Archimedean approach. Though the Ding functional is continuous on $\mathcal E^1$
with respect to $d_1$, the Mabuchi functional is only lower semi-continuous (c.f. Lemma 4.3 in \cite{BBEGZ19}).
Darvas-Zhang replace the entropy term $H(u)$ in the Mabuchi energy by
$$ H^{\beta}(u) = \sup_{v \in \mathcal H_\omega} \left( - \log \int_X e^{\beta(v-u)}\frac{\omega^n}V + \beta(E(v) - E(u))\right),$$
and define $M^\beta : \mathcal H_\omega \to \mathbf R$ by
$$ M^\beta(u) = H^\beta(u) - n E_{Ric}(u) + \bar S E(u).$$
It is shown that $H^\beta(u) \le H(u)$ and that $H^\beta(u) \nearrow H(u)$ as $\beta \nearrow \infty$, and thus
$M^\beta(u) \nearrow M(u)$ as $\beta \nearrow \infty$. $M^\beta$ is continuous with respect to $d_1$. The properness
of $M$ is equivalent to the properness of $M^\beta$ for $\beta$ big enough.

To an ample test configuration $(\mathcal X, \mathcal L)$, we can
associate a $C^{1,1}$-geodesic ray $\{u_t\}_{t \in [0,\infty)} \subset \mathcal H_\omega^{1,1}$
(\cite{PH07}, \cite{PH10}, \cite{RWN14}, \cite{Berman16}, \cite{CTW18}), and introduce an 
invariant 
$$ M^\beta(\mathcal X, \mathcal L) := \liminf_{t \to \infty} \frac{M^\beta(u_t)}t.$$
Given $\beta > 0$, we say $(X, L)$ is uniformly $K^\beta$-stable if there exists $\epsilon > 0$ such that
$$ M^\beta(\mathcal X, \mathcal L) \ge \epsilon J_{na}(\mathcal X, \mathcal L).$$
The following is the first characterization only using ample test configurations as test objects.
\begin{theorem}[Darvas-Zhang \cite{DarvasZhang25}]
Let $(X,L)$ be a polarized manifold. The following are equivalent:\\
(i) $(X, L)$ admits a unique cscK metric in $c_1(L)$;\\
(ii) $(X, L)$ is uniformly $K^\beta$-stable for some $\beta > 0$.
\end{theorem}
Darvas-Zhang also give an algebraic description of $M^\beta(\mathcal X, \mathcal L)$, and an interpretations in terms of models.
	


\vspace{0.3cm}

Added after submission : After submission of this survey paper, a paper by Trusiani \cite{Trusiani26} appeared. It is proved in this paper that a polarized smooth projective variety $(X,L)$ admits a cscK metric if and only if it is $\mathrm{Aut}_0(X,L)$-uniformly K-stable.

\section{Examples of Conformally K\"ahler, Einstein-Maxwell metrics}
Weighted extremal metrics contain various special metrics 
and each has its own unique feature.
In this section, we review conformally K\"ahler, Einstein-Maxwell metrics. 

The following are known examples of conformally K\"ahler, Einstein metrics.

- Page \cite{Page78} : A conformally K\"ahler, Einstein metric on the one-point-blow-up of $\bfC\bfP^2$.

- Chen-LeBrun-Weber \cite{ChenLeBrunWeber} : A conformally K\"ahler, Einstein metric on the two-points-blow-up of $\bfC\bfP^2$.

-  Apostolov-Calderbank-Gauduchon \cite{ACG15} : Conformally K\"ahler, Einstein metrics on 4-orbifolds.

- B\'erard-Bergery \cite{BB82} : Conformally K\"ahler, Einstein metrics on $\mathbf P^1$-bundles over Fano K\"ahler-Einstein manifolds.

Next let us review known examples of non-Einstein cKEM metrics.

- LeBrun \cite{LeBrun16} : Ambitoric examples on $\mathbf C\mathbf P^1 \times 
\mathbf C\mathbf P^1$ and the one-point-blow-up
of $\mathbf C\mathbf P^2$.

- Koca-T{\o}nnesen-Friedman \cite{KT} : Examples on ruled surfaces of higher genus.

- Viza de Souza \cite{Viza de Souza21} : A more example on the one-point-blow-up
of $\mathbf C\mathbf P^2$.

The last example by Viza de Souza is obtained by a principle of volume minimization.
Let us review on the volume minimization. In \cite{MSY2}, 
Martelli-Sparks-Yau formulated the volume minimization to obtain
Sasaki-Einstein metrics.
They considered a volume functional on the space of toric Sasakian structures. This
functional is convex and proper, and thus there is a critical point, which corresponds 
to a Sasakian structure with vanishing Sasaki-Futaki invariant. It is shown
in \cite{FOW} that toric Sasaki manifold with vanishing Sasaki-Futaki invariant has a
Sasaki-Einstein metric. This way, we can find many examples of Sasaki-Einstein metrics.

For K\"ahler-Ricci solitons, there has been a similar story due to Tian-Zhu \cite{TianZhu02}.
They consider a volume functional on the space of holomorphic Killing vector fields, which is convex and proper.
The critical points correspond to a Killing vector field for which an obstruction to the
existence of K\"ahler-Ricci solition vanishes.

Motivated by these works, Ono and the author considered in \cite{FO17}
a volume functional on the space of holomorphic Killing vector fields
whose critical points correspond to the conformal factor $f$ in section 3.
The volume functional is neither convex nor proper in this case. 
For example, on the one-point-blow-up of $\mathbf C\mathbf P^2$, there is a K\"ahler class for which there 
are two critical points. One of them correspond to 
the examples of LeBrun \cite{LeBrun16} and the other had not been known.
Viza de Souza succeeded our study and showed that the other one certainly corresponds
to a cKEM metric. This is how Viza de Souza's example was found.

See \cite{FO_ICCM_Notices19}, \cite{ACL21}, \cite{AJL21}, \cite{futaki26JDGSurvey} for more on weighted extremal metrics and the special case weighted solitons.

\end{document}